\documentclass{mcom-l}

\usepackage{amsmath,amsfonts,amssymb,amsthm}
\usepackage{graphicx,multirow,rotate,color,subfigure}
\usepackage{lipsum,verbatim}

\newcommand\Tstrut{\rule{0pt}{2.6ex}}         

\usepackage{color,colortbl}
\definecolor{red}{rgb}{1,0,0}
\definecolor{green}{rgb}{0,1,0}
\definecolor{blue}{rgb}{0,0,1}
\definecolor{magenta}{rgb}{1,0,1}

\definecolor{ciano}{rgb}{0,1,1}
\definecolor{magenta}{rgb}{1,0,1}
\definecolor{amarelo}{rgb}{1,1,0}
\definecolor{iw}{rgb}{0.7,0.93,0.36}  
\newcommand{\iw}{\cellcolor{iw}}
\newcommand{\ia}{\cellcolor{amarelo}}

\begin{document}

\title[Gauss periods are minimal discriminant polynomials]{Gauss periods are minimal polynomials for totally real cyclic fields of prime degree}

\author{Jason A.C.~Gallas}
\address{Max-Planck Institut PKS, N\"othnitzer Str. 38, 01187 Dresden, Germany}
\email{jgallas@pks.mpg.de}
\address{Dept.~of Appl.~Math., Complexity Sciences Center, 9225 Collins Avenue
         Suite 1208,  Surfside, FL 33154, USA}
\address{Instituto de Altos Estudos da Para\'\i ba, 
          Rua Silvino Lopes 419-2502\\ 
          58039-190 Jo\~ao Pessoa, Brazil}


\begin{abstract}
We report extensive computational evidence that Gauss period equations
are  minimal discriminant polynomials for primitive elements representing
Abelian (cyclic) polynomials of prime degrees $p$.
By computing 200 period equations up to $p=97$, we
significantly extend tables in the  compendious number fields database
of Kl\"uners and Malle. Up to $p=7$, period equations reproduce known
results  proved to have minimum discriminant.
For $11\leq p\leq 23$, period equations coincide with 53 known but
unproved cases of minimum discriminant in the database, and fill a gap
of 19 missing cases.
For $29\leq p\leq 97$, we report 128 not previously known cases, 16 of them
conjectured to be minimum discriminant polynomials of Galois group $pT1$.
The significant advantage of period equations is that they all may be obtained
analytically using a procedure that works for fields of arbitrary degrees,
and which are extremely hard to detect by systematic numerical search.
\end{abstract}

\maketitle 

\section{Introduction}

The determination of number field representations by a minimal polynomial
of primitive elements of interest keeps drawing attention as may be seen from
Refs.~\cite{an76}--\cite{ba19}, and many additional references given therein.
Such polynomials are notoriously difficult to determine, particularly for
degrees higher than 10, say.
For the 9 primes $p\leq23$, they are available in the
compendious and user-friendly tables of Kl\"uners and Malle, 
containing more than $5\times10^6$ polynomials, and
obtained by extensive systematic computer search \cite{km}--\cite{m06},

The aim of this paper is to report the computation of number field
representations of minimal polynomial for the 25 primes $p < 100$.
In contrast with previous works, such polynomials are obtained here totally
analytically. As shown below, this was possible after realizing that 
minimal polynomials coincide with the {\sl period equations} described
in \S 343 of the  {\sl Disquisitiones Arithmetic\ae} \cite{da}:
{\it Omnes radices $\Omega$ in certas classes (periodos) distribuuntur}.

As an extra bonus, period equations allow minimal polynomials
to be computed {\sl analytically for arbitrarily large primes $p$},
well beyond what is possible to probe with systematic numerical search,
due to the very large number and magnitudes of the polynomial
coefficients involved.
Some representative large coefficients are given in Table \ref{tab:tab04} below.
Here we focus on totally real Abelian (cyclic) fields because of their
relevance as equations of motion for certain classes of dynamical systems,
namely in the so-called partition generating limit of discrete-time
quadratic and quartic dynamics \cite{carriers,wolf}.

\section{Analytical expressions for period equations and field discriminants}

Consider the set $\Omega$ of the $p-1$ complex roots of the equation $x^p-1=0$.
The key to obtain Gaussian period equations is to distribute the $p-1$ roots
in $\Omega$ into certain sums called  {\sl periods} \cite{da,bach}.
This is done by choosing two positive rational integers $e$ and $f$
such that $ef=p-1$,
and partitioning the roots in $\Omega$ into $e$ disjoint classes, thereby
forming $e$ periods $\eta_i$, each one consisting of a sum of $f$ roots.

\setlength{\tabcolsep}{7pt} 
\begin{table}[!tbh]
\caption[]{Field discriminants for families of totally real cyclic fields
of prime degree $p$ and Galois group $pT1$, computed with Eq.~(\ref{disco}).
Non-highlighted discriminants coincide with values in the tables of  
Kl\"uners and Malle \cite{km}, while highlighted discriminants are missing
in their tables.
Minimal polynomials representing the 19 new fields are given in
Table \ref{tab:tab02}.
Labels $\ell$ refer to the order in which polynomials are generated when
$f$ is varied in $ef=p-1$ along rows.
The $\ell=1$ column lists proven and unproven {\sl minimum} field discriminants.
}
 \label{tab:tab01}
{\small
\begin{tabular}{c cccccc cccc}
\hline
$p$ & $\ell=1$ & $\ell=2$ & $\ell=3$ & $\ell=4$ & $\ell=5$ & $\ell=6$ & $\ell=7$ & $\ell=8$\\
\hline
2 & $5$& $13$& $17$& $29$& $37$& $41$& $53$ & $61$ \\
3 & $7^2$&  $13^2$& $19^2$& $31^2$& $37^2$& $43^2$& $61^2$ & $67^2$ \\
5 & $11^4$  & $31^4$ & $41^4$ & $61^4$ & $71^4$ & $101^4$ & $131^4$ & $151^4$\\ 
7 & $29^6$ & $43^6$ & $71^6$ & $113^6$ & $127^6$ & $197^6$ & $211^6$
           & $239^6$\\ 
11& $23^{10}$ & $67^{10}$ & $89^{10}$ & \iw{$199^{10}$} & \iw{$331^{10}$} &
           \iw{$353^{10}$} & \iw{$397^{10}$} & \iw{$419^{10}$}\\
13& $53^{12}$ & $79^{12}$ & \iw{$131^{12}$}& \iw{$157^{12}$}&
     \iw{$313^{12}$}& \iw{$443^{12}$}& \iw{$521^{12}$} & \iw{$547^{12}$}\\
17& $103^{16}$& $137^{16}$& $239^{16}$& \iw{$307^{16}$}& \iw{$409^{16}$}
             &\iw{$443^{16}$}& \iw{$613^{16}$} &  \iw{$647^{16}$} \\
19& $191^{18}$& $229^{18}$& $419^{18}$& $457^{18}$ & $571^{18}$& \iw{$647^{18}$}&
        \iw{$761^{18}$} & \iw{$1103^{18}$}\\
23& $47^{22}$& $139^{22}$& $277^{22}$& $461^{22}$& $599^{22}$& $691^{22}$&
          $829^{22}$ & $967^{22}$\\
\hline
\end{tabular}
}
\end{table}

Let $r$ be any complex root in $\Omega$, say $r=e^{2\pi i/p}$,
and $g$ a primitive root of $p$. 
Then, the periods  $\eta_i$ are given by the sums \cite{da,bach,reu}
\begin{equation}
  \eta_i = \sum_{k=0}^{f-1} r^{g^{ke+i}}, \qquad i=0,1,\cdots,e-1,
\end{equation}
namely, $e$ sums of $f$ distinct complex roots suitably selected from $\Omega$ 
\begin{alignat*}{8}
  \eta_0 &=& r &+&r^{g^e} &+&r^{g^{2e}} &+&\ \cdots\  &+&r^{g^{(f-1)e}},\\
  \eta_1 &=& r^{g} &+&r^{g^{e+1}} &+& r^{g^{2e+1}} &+&\ \cdots\ &+&\ r^{g^{(f-1)e+1}},\\
  &\; \; \vdots&\\
  \eta_{e-1} &=&\ r^{g^{e-1}} &+&\ r^{g^{2me1}} &+&\ r^{g^{3e-1}} &+&\ \cdots\
                  &+&\ r^{g^{fe-1}}.
\end{alignat*}  

For a fixed $e$, the {\sl period equation} $\psi_e(x)$ of degree $e$
is defined as \cite{da,bach,reu}
\begin{equation}
  \psi_e(x) = \prod_{k=0}^{e-1} (x-\eta_{k})
            = x^e+x^{e-1}+\alpha_2x^{e-2}+\cdots+ \alpha_e,    \label{theta}
\end{equation}
where the $\alpha_\ell$ turn out to be integers.
The corresponding field discriminant $\Delta_e$ of $\psi_e(x)$ is \cite{jg20}
\begin{equation}\label{disco}
    \Delta_e =
    \begin{cases}
  -p^{e-1}, & \hbox{if } (e-1) \  \hbox{mod } 4 =1 {\ \ \ \rm and\ \ \ }
                       f {\ \rm mod\ } 2 = 1, \\
  \phantom{-}p^{e-1}, & \hbox{if\ \ otherwise}.
    \end{cases}
\end{equation}
Equivalently, Eq.~(\ref{disco}) may be also written as
\begin{equation}\label{disco2}
    \Delta_e =
    \begin{cases}
 (-1)^{n_P} p^{e-1}, & \hbox{if} \ (e-1) \ \hbox{mod } 4 =1, \\
      p^{e-1},           & \hbox{if\ \ otherwise},
    \end{cases}
\end{equation}
where $n_P$ is the number of {\it pairs} of complex roots of $\psi_e(x)$.

\setlength{\tabcolsep}{7pt} 
\begin{table}[!tbh]
\caption[]{Polynomials of Galois group $pT1$ derived analytically which
fill gaps in tables of Kl\"uners and Malle \cite{km}.
For reference, the highlighted polynomials in each group may be compared
with analogous ones in the tables Kl\"uners and Malle.
$\Delta_e$ is the field discriminant.
}
 \label{tab:tab02}
{\tiny
  \begin{tabular}{c l}
\hline  
$p,\Delta_e$  & Polynomial coefficients ordered by decreasing powers of $x$\\
\hline 
\ia{$11,89^{10}$}&\ia $ 1, 1, -40, -19, 482, 84, -2185, 102, 3152, -781, 57, -1$\Tstrut\\
$11,199^{10}$ & $1, 1, -90, -115, 2349, 943, -26327, 21284, 102168, -217794, 148930, -30647$\\
$11,331^{10}$ & $1, 1, -150, -402, 6577, 28617, -62124, -475464, -343344, 1913488, 4015168, 2287616$\\
$11,353^{10}$ & $1, 1, -160, -525, 6066, 26034, -48369, -265374, -42966, 405001, 63189, -170569$\\
$11,397^{10}$ & $1, 1, -180, -13, 11655, -12159, -316973, 720142, 2670510, -10551746, 10752776, -3098903$\\
$11,419^{10}$    &$1, 1, -190, -547, 10985, 51221, -141765, -1224028, -2399676, -1263744, 873500, 785489$\\
\hline 
\ia{$13,79^{12}$}&\ia$1, 1, -36, -77, 365, 1193, -617, -5541, -4414, 4575, 6321, -411, -2196, -293$\Tstrut\\
$13,131^{12}$&  $1, 1, -60, -27, 1199, 33, -9610, 3352, 33548, -20328, -47723, 34869, 21271, -15667  $\\
$13,157^{12}$&  $1, 1, -72, -129, 1672, 3386, -16810, -32367, 81708, 121902, -196272, -127412, 217458, -61399  $\\
$13,313^{12}$&  $1, 1, -144, -161, 6530, 9620, -109398, -196143, 512628, 917970, -650724, -1134730, 253950, 409375  $\\
$13,443^{12}$&  $1, 1, -204, 181, 10752, -9116, -208418, 161679, 1686466, -1207646, -4904338, 3051848, 896956, -144209  $\\
$13,521^{12}$&  $1, 1, -240, 293, 19153, -45777, -616830, 1795569, 7791196, -23224049, -29107980, 68466088, 31673025, -4516075  $\\
$13,547^{12}$&  $1, 1, -252, -1123, 15626, 107844, -204415, -3094114, -4853400, 22393129, 91453411, 116380476, 47088126, 1165671  $\\
\hline  
\ia{$17,239^{16}$}&\ia$1, 1, -112, -47, 3976, 4314, -64388, -136247, 422013, 1631073, 411840, -5840196, -11894369, -10635750$,\Tstrut\\
 &\qquad$-4739804,-938485, -54850, -619  $\\
$17,307^{16}$&  $1, 1, -144, -241, 6894, 14938, -127923, -323969, 847982, 2194186, -2617873, -6091397, 3745755, 7069429, -1600190$,\\
   &\qquad$-3100257, -220118, 208777  $\\
$17,409^{16}$&  $1, 1, -192, -273, 14752, 28028, -571107, -1411675, 11275657, 36814399, -91832077, -461179352, -109192148$,\\
 & \qquad$\phantom{-}1929139488,3679722325, 2767754010, 828153361, 45886883$\\
$17,443^{16}$&  $1, 1, -208, 17, 15287, -13881, -487578, 703261, 6754359, -10540902, -41136753, 57683825, 92010954, -95287840$,\\
 &\qquad$-17501435, 25026156, -563260, -1246103$\\
$17,613^{16}$&  $1, 1, -288, -265, 26034, 40228, -875968, -2022008, 8464009, 27681440, -8855367, -101412811, -87313302, 38624139,$\\
&  \qquad$\phantom{-}67164168,7149746, -7878215, 664471$\\
$17,647^{16}$&  $1, 1, -304, -1117, 25631, 126439, -773932, -4360454, 10731832, 64676368, -79260104, -441919082, 345306489$,\\
&\qquad\phantom{\ -}$1259087517,-718017711, -1025767171, 183044979, 202031659$\\
\hline  
\ia{$19,571^{18}$}&\ia$1, 1, -270, -351, 28987, 41181, -1648168, -2453428, 55186847, 85340779, -1133336624, -1826548777, 14258200659$,\Tstrut\\
&\qquad\phantom{\ -}$24108274876,-104945874488, -187074906809, 398834944226, 749021546949, -549551190705, -1072348621073$\\

$19,647^{18}$&  $1, 1, -306, 79, 29370, -130, -1361758, -704483, 34401066, 34840236, -488507308, -678792838, 3902105097, 6509889712$,\\
&\qquad$-16762277726, -31790175776, 33454527221, 71781895812, -18736713884, -50781641759  $\\

$19,761^{18}$&  $1, 1, -360, 173, 45376, -58762, -2558302, 4227138, 70534890, -120121397, -973700212, 1501612590, 6678374954$,\\
&\qquad$-8595059019,-21259099080, 21796436285, 27241052007, -18814754704, -12659238391, 3483379661$\\
$19,1103^{18}$&  $1, 1, -522, -504, 101486, 128414, -9620928, -15777284, 481094177, 987900729, -12621266990, -31913538892, 160258304904$,\\
&\qquad\phantom{\ -}$502132433072, -773726987936, -3469327268928, -293071915904, 7722407112704, 6083919470592, -276663107584  $\\
\hline
\end{tabular}
}
\end{table}

The explicit calculation of period equations presents no more difficulty
than the inevitable length.
Analytically, the method works fine for low-degree equations.
For instance,
it is a standard textbook example to solve $x^{17}-1=0$ by radicals.
But the method collapses in a tangle of calculations as soon as it is applied
to higher degrees, leading to no convenient formulas to solve polynomials
by radicals.
Fortunately, the method is not difficult to implement in computer algebra
systems.
Using the equations above, we wrote a procedure in MAPLE 18 (X86 64 Linux),
which was used to produce the polynomials and results summarized in our
tables below.

\section{Results}

Table \ref{tab:tab01} summarizes  field discriminants calculated for the
$\psi_e(x)$ polynomials for $p\leq23$. The non-highlighted discriminants
arise from polynomials that coincide with the ones included in the
Tables of Kl\"uners and Malle \cite{km}.
The 19 highlighted discriminants are from polynomials missing in their
tables.
In Table \ref{tab:tab01}, the labels $\ell$ refer to the order in which
polynomials are generated when $f$ is varied in $ef=p-1$ along rows.

Table \ref{tab:tab02} reports polynomial coefficients for the 19 highlighted
cases in Table \ref{tab:tab01}.
Coefficients are arranged in order of decreasing powers of $x$, for
instance, the topmost polynomial for $p=11$ is
\[  {x}^{11}+{x}^{10}-40\,{x}^{9}-19\,{x}^{8}+482\,{x}^{7}+84\,{x}^{6}
           -2185\,{x}^{5}+102\,{x}^{4}+3152\,{x}^{3}-781\,{x}^{2}+57\,x-1. \]

Table \ref{tab:tab03}  summarizes discriminants computed for 128 new
Gaussian period equations $\psi_e(x)$.
Note the relatively large magnitudes of the discriminants,
e.g.~$43651^{96}$ is a number of  $446$ digits.
Of particular interest are the discriminants listed in the $\ell=1$ column
which, based on what happens in all previous cases, are conjectured to be
minimum field discriminant polynomials. These polynomials are given
explicitly in Table \ref{tab:tab04}.
The very large magnitude of their coefficients should make them rather
difficult, if not impossible, to determine by systematic computer searches.
Note that polynomials for additional primes $p$ are not difficult to be
determined, thanks to the analytical expression for $\psi_e(x)$.

In Tables \ref{tab:tab02} and \ref{tab:tab04}, the $\ell=1$ period equations
are also exceptional in that they are {\sl monogenic} polynomials,
meaning that their ordinary discriminant $\Delta_\psi$
of $\psi_e(x)$ coincides with their field discriminant $\Delta_e$,
a rather rare occurrence \cite{gaal} that misled Dedekind for some time
\cite{dedek}.
For other values of $\ell$ in Tables \ref{tab:tab01} and \ref{tab:tab03} one has
\[ \Delta_\psi = k^2\;\Delta_e, \qquad k^2\neq 1. \]

\setlength{\tabcolsep}{7pt} 
\begin{table}[!tbh]
\caption[]{Field discriminants for 128 new $\psi_e(x)$
  period equations  
  which significantly extend the known results in Table \ref{tab:tab01}.
Labels $\ell$ refer to the order in which polynomials are generated when
$f$ is varied in $ef=p-1$ along rows.
Discriminants on the $\ell=1$ column are conjectured to be
{\sl minimum} field discriminants. The corresponding
   polynomials are given explicitly in  Table \ref{tab:tab04}.}  
 \label{tab:tab03}
{\small
\begin{tabular}{c cccccc cccc}
\hline
$p$ & $\ell=1$ & $\ell=2$ & $\ell=3$ & $\ell=4$ & $\ell=5$ & $\ell=6$ & $\ell=7$ & $\ell=8$\\  
\hline
29& \ia$59^{28}$& $233^{28}$& $349^{28}$& $523^{28}$
                   & $929^{28}$& $1103^{28}$& $1277^{28}$& $1451^{28}$ \\
31 & \ia$311^{30}$& $373^{30}$& $683^{30}$& $1117^{30}$& $1303^{30}$
   & $1427^{30}$& $1489^{30}$& $1613^{30}$ \\
37 & \ia$149^{36}$& $223^{36}$& $593^{36}$ & $1259^{36}$
   & $1481^{36}$& $1777^{36}$& $1999^{36}$& $2221^{36}$\\
41 & \ia$83^{40}$& $739^{40}$& $821^{36}$& $1231^{40}$
   & $1559^{40}$& $1723^{40}$& $2297^{40}$& $2543^{40}$\\
43 & \ia$173^{42}$& $431^{42}$& $947^{42}$& $1033^{42}$
   & $1291^{42}$& $1549^{42}$& $1721^{42}$& $1979^{42}$\\
47 & \ia$283^{46}$& $659^{46}$& $941^{46}$& $1129^{46}$
   & $1223^{46}$& $1693^{46}$& $1787^{46}$& $2069^{46}$\\
53 & \ia$107^{52}$& $743^{52}$& $1061^{52}$& $1697^{52}$
   & $2333^{52}$& $2969^{56}$& $3181^{52}$& $3499^{52}$\\
59 & \ia$709^{58}$& $827^{58}$& $1063^{58}$& $1181^{58}$
   & $1889^{58}$& $2243^{58}$& $2833^{58}$& $3187^{58}$\\
61 & \ia$367^{60}$& $977^{60}$& $1709^{60}$& $1831^{60}$
   & $2441^{60}$& $3539^{60}$& $4027^{60}$& $4271^{60}$\\
67 & \ia$269^{66}$& $1609^{66}$& $1877^{66}$& $2011^{66}$
   & $3083^{66}$& $3217^{66}$& $4021^{66}$& $4289^{66}$\\
71 & \ia$569^{70}$& $1279^{70}$& $2131^{70}$& $2273^{70}$
   & $2557^{70}$& $2699^{70}$& $4261^{70}$& $5113^{70}$\\
73 & \ia$293^{72}$& $1607^{72}$& $1753^{72}$&  $3359^{72}$& $3797^{72}$
   & $3943^{72}$& $4673^{72}$& $5987^{72}$\\  
79 & \ia$317^{78}$& $2213^{78}$& $2687^{78}$& $3319^{78}$
   & $3793^{78}$& $5531^{78}$& $6163^{78}$& $6637^{78}$ \\
83 & \ia$167^{82}$&  $1163^{82}$& $1993^{82}$& $2657^{82}$
   & $4483^{82}$& $5147^{82}$& $5479^{82}$& $6143^{82}$\\
89 & \ia$179^{88}$& $4273^{88}$& $5519^{88}$& $6053^{88}$& $7477^{88}$
   & $8011^{88}$& $9257^{88}$& $9781^{88}$\\
97 & \ia$389^{96}$&  $3881^{96}$& $10477^{96}$& $11059^{96}$
   & $23087^{96}$& $25997^{96}$& $42293^{96}$& $43651^{96}$  \\
\hline
\end{tabular}
}
\end{table}

\section{Conclusions}

The main message of this paper is the identification of minimal polynomials
for totally real cyclic fields as being Gaussian period equations.
This fact allowed the easy analytical determination of
period equations and field discriminants for the 25 cyclic fields of
prime degree $p\leq100$. Such period equation were found to
coincide with all known minimum discriminant polynomials for totally real
cyclic fields  and to significantly extend the minimal polynomials by 128
new cases.
Up to $p=7$, the 4 period equations coincide with known polynomials proven
to have minimum discriminants.
For $11\leq p\leq 23$, the 5 equations coincide with tabulated polynomials
having known but still unproven minimum discriminant.
The remaining 16 equations,
for $29\leq p\leq 97$, we reported 128 not previously known minimal polynomials,
16 of them conjectured to also have minimum discriminants
The findings reported here suggest that polynomials representing totally
real number fields with minimum discriminants are given by Gaussian period
equations which may be determined to very high degrees with no major difficulty.
Although we limited our tables to just $\ell=8$ cases, we have investigated
much higher  $\ell$ values with no difficulty, always corroborating the
results described above.

The ability of the Gaussian period equations to reproduce minimal polynomials
is almost uncanny,
and we hope this work to set the agenda for the study of their general
properties and new applications.
The identification of Gaussian period equations as representatives of minimal
discriminant fields was also found to be valid for fields of arbitrary degrees,
and for fields of signatures other than totally real.
These results will be presented elsewhere.

\section*{Acknowledgments}

This work was partially supported by the Max-Planck Institut PKS, Dresden,
Germany, and
CNPq, Brazil,  Grant No.~PQ-305305/2020-4.


\appendix
\section{Table \ref{tab:tab04}: Gaussian period equations extending  Table \ref{tab:tab01}}

\setlength{\tabcolsep}{7pt} 
\begin{table}[!tbh]
\caption[]{Analytically derived coefficients of Gaussian period equations
  which significantly extend the known cases shown in Table \ref{tab:tab01}
of totally real cyclic fields of prime degree and Galois group $pT1$.
These polynomials are conjectured to have minimum field discriminants
$\Delta_e$. Numbers in parenthesis give the number of digits of $\Delta_e$.
}
 \label{tab:tab04}
{\tiny
  \begin{tabular}{c l}
\hline  
$p,\Delta_e$  & Polynomial coefficients ordered by decreasing powers of $x$\\
\hline
$29,59^{28}$ & $1, 1, -28, -27, 351, 325, -2600, -2300, 12650, 10626, -42504, -33649, 100947, 74613, -170544, -116280, 203490, 125970$,\Tstrut\\
(50)&\qquad$-167960, -92378, 92378, 43758, -31824, -12376, 6188, 1820, -560, -105, 15, 1  $\\
$31,311^{30}$&  $1, 1, -150, -117, 9434, 4958, -328880, -78328, 7098810, -507925, -100297691, 39217806, 950942678, -677567200, -6056859921$,\Tstrut\\
(75)&\qquad$\phantom{-}6293114670, 25303767350, -35315191500, -65145555320, 122001770825, 87706314807$,\\
&\qquad$-253609179978, -19776440614, 299693882286, -98200795275, -175359843704, 113820457794, 30260338943$,\\
&\qquad$-38729335872, 5811762482, 1056847214, 26458109$\\
$37,149^{36}$&  $1, 1, -72, -67, 2278, 1964, -41849, -33262, 497142, 362207, -4027127, -2672215, 22871279, 13719683, -92273266, -49606380$,\\
(79)&\qquad$\phantom{-}265291071, 126402144, -540999671, -224635365, 773506030, 271832880, -761221384, -214848879, 502240682, 103923082$,\\
&\qquad$-213804070, -27920639, 55105551, 3603129, -7708785, -235974, 497852, 22640, -12819, -1020, 28, 1$\\
$41,83^{40}$&  $1, 1, -40, -39, 741, 703, -8436, -7770, 66045, 58905, -376992, -324632, 1623160, 1344904, -5379616, -4272048, 13884156$,\\
(77)&\qquad$\phantom{-}10518300, -28048800, -20160075, 44352165, 30045015, -54627300, -34597290, 51895935, 30421755, -37442160, -20058300$,\\
&\qquad$\phantom{-}20058300, 9657700, -7726160, -3268760, 2042975, 735471, -346104, -100947, 33649, 7315, -1540, -210, 21, 1$\\
$43,173^{42}$&  $1, 1, -84, -79, 3160, 2786, -70521, -58076, 1042773, 798942, -10810701, -7671648, 81126975, 53056499, -448758019, -268953093$,\\
(94)&\qquad$\phantom{-}1846875875, 1007873658, -5671315301, -2797394685, 12962901258, 5730420663, -21895326590, -8590544711, 27001558938$,\\
&\qquad$\phantom{-}9297003415, -23896748020, -7121607714, 14834334408, 3755010451, -6269987218, -1310830451, 1732796449, 287529242$,\\
&\qquad$-294221159, -37041429, 27505221, 2599666, -1134285, -90620, 12893, 370, -54, 1$\\
$47,283^{46}$&  $1, 1, -138, -315, 8338, 29804, -276833, -1433626, 5033859, 41190458, -30657314, -748097961, -742659788, 8506344013, 21519259357$,\\
(113)&\qquad$-52948548811, -268879641855, 30332528938, 1920252236103, 2430736233424, -7367030656288, -21401598866455, 5373046913681$,\\
&\qquad$\phantom{-}89240242581627, 85735098102709, -174332690558567, -418760640969237, -24016008538438, 845143941649693, 850093833789498$,\\
&\qquad$-563502754038610, -1652337279635119, -688139958495907, 1132046847057208, 1399925903803443, 129366579867529$,\\
&\qquad$-739733924950208, -464611724348759, 61382685657965, 169596960618209, 46703998783537, -17972540150505, -11555848466611$,\\
&\qquad$-591153994800, 757301830397, 134933083169, -1756074840, -132954859$\\
$53,107^{52}$&  $1, 1, -52, -51, 1275, 1225, -19600, -18424, 211876, 194580, -1712304, -1533939, 10737573, 9366819, -53524680, -45379620, 215553195$,\\
(106)&\qquad$\phantom{-}177232627, -708930508, -563921995, 1917334783, 1471442973, -4280561376, -3159461968, 7898654920, 5586853480$,\\
&\qquad$-12033222880, -8122425444, 15084504396, 9669554100, -15471286560, -9364199760, 12875774670, 7307872110, -8597496600$,\\
&\qquad$-4537567650, 4537567650, 2203961430, -1855967520, -818809200, 573166440, 225792840, -129024480, -44352165, 20160075$,\\
&\qquad$\phantom{-}5852925, -2035800, -475020, 118755, 20475, -3276, -351, 27, 1$\\
$59,709^{58}$&  $1, 1, -348, -773, 54710, 181658, -5108670, -22380037, 313864889, 1735112733, -13244902303, -92489323874, 384874485628$,\\
(166)&\qquad$\phantom{-}3563558672105, -7224424298771, -102284588313436, 59507208676899, 2226843428814170, 1091206962010061$,\\
&\qquad$-37110651573131926, -50769946526736009, 473955610044577175, 1034887463549316633, -4599282007018108614$,\\
&\qquad$-14157628456571325637, 33064370721903906102, 142351088507771542259, -164489546968394181788, -1088242867211629932856$,\\
&\qquad$\phantom{-}434724808533168556934, 6421863971239504423218, 821416241786189985129, -29458267038177484817142$,\\
&\qquad$-15268230232813063311459, 105361097493218843884743, 86357419790172643990359, -294180307521278014201331$,\\
&\qquad$-310904964780280912527690, 641819365276550155193989, 795098801640263156488019, -1096213896340607733888156$,\\
&\qquad$-1490164806719968998511703, 1470729980576933657823174, 2055743062031065158143660, -1554777323587753672707889$,\\
&\qquad$-2062294848329571542288479, 1289894730680748454358450, 1461538845678796120527548, -820373077133361026883291$,\\
&\qquad$-693872641742794094654444, 378259269617360523640282, 199906595981055748193122, -114544083882472863320184$,\\
&\qquad$-27852340198822864548711, 19330914962494272607266, 388755699390176021351, -1296740570702255596842$,\\
&\qquad$\phantom{-}197731971687386122768, -8442821363163035954, -4584804110356861$\\
$61,367^{60}$&  $1, 1, -180, -413, 14645, 52673, -688601, -3528285, 19601342, 147066683, -294411763, -4094132898, -515253399, 77730183943$,\\
(154)&\qquad$\phantom{-}136127326190, -976171853534, -3440181895427, 6820700430038, 49495243759604, 5003873089663, -454773118963937$,\\
&\qquad$-708332454647506, 2503079468369386, 8524707843825128, -4353506873916376, -55965866407852931, -50586504885325349$,\\
&\qquad$\phantom{-}208588686855790630, 488766704388771586, -248616968576483612, -2111939395471861259, -1621839831538649260$,\\
&\qquad$\phantom{-}4656502567948889185, 9658116194858556310, -1497012549402345042, -23738135773030440428, -20637249171326524602$,\\
&\qquad$\phantom{-}24569280463351126387, 57740462789299170679, 14579305016676768566, -65771381500668227611, -72094504453297506055$,\\
&\qquad$\phantom{-}13570829828275968614, 78797457647783506537, 45464576927398898131, -26397733136489030972, -45803590606881784425$,\\
&\qquad$-13479501565046232768, 13463013755243374867, 12223486303700829006, 1531995617799024854, -2479700985326385279$,\\
&\qquad$-1203853517042779781, 45547818199089944, 173724631575165523, 33343190173223859, -8176066599449700$,\\
&\qquad$-3216273109307699, 787254613563, 91700724228752, 3245424268762, -934868803397$\\
$67,269^{66}$&  $1, 1, -132, -127, 8128, 7514, -310489, -275412, 8255097, 7015522, -162469429, -132074732, 2457934103, 1908054027, -29302580527$,\\
(161)&\qquad$-21682477578, 280016358684, 197116015784, -2170586468374, -1450574505983, 13762204743256, 8711498507598, -71769908578626$,\\
&\qquad$-42926578989207, 308917467707281, 174121538109300, -1099328139842573, -582252857213765, 3235005088514751$,\\
&\qquad$\phantom{-}1605012951198723, -7862532891652511, -3641973093600359, 15743247597494598, 6784769007620029, -25870357696878572$,\\
&\qquad$-10338065796135117, 34707160257704653, 12822873516191615, -37759031607934407, -12874525720258242, 33032386678911041$,\\
&\qquad$\phantom{-}10396393258198832, -22995232909296191, -6702394382020432, 12575711941082825, 3418837407991456, -5318199438079755$,\\
&\qquad$-1363375197739315, 1705561117819514, 417633543753729, -404848315593940, -95674661443562, 68987247039926$,\\
&\qquad$\phantom{-}15749316188706, -8123019166449, -1758233811824, 632045599194, 122538815988, -31180078280, -4726423766$,\\
&\qquad$\phantom{-}955488656, 84545195, -16930746, -385257, 133369, -3611, -84, 1$\\
$71,569^{70}$&  $1, 1, -280, -259, 36362, 31280, -2918219, -2341868, 162684328, 121979735, -6712586533, -4705153431, 213256136347, 139745967919$,\\
(193)&\qquad$-5356937166534, -3281290140692, 108406810097566, 62068478295763, -1791341741101037, -959090929705315, 24407831840018355$,\\
&\qquad$\phantom{-}12233821563475166, -276147939713855185, -129846912904598319, 2606654089690226864, 1153616129990347133$,\\
&\qquad$-20587950147002420775, -8616230414933743750, 136232462096831704169, 54248999168315159885, -755063095484716789902$,\\
&\qquad$-288262405489471572571, 3499387584477151458385, 1291958342775348433576, -13519262273839572481549$,\\
&\qquad$-4871234681829364379697, 43337921624309469265580, 15374832177489501594722, -114567325220911193557243$,\\
&\qquad$-40314399744765933113883, 247818220129501437242168, 86899209178511442518249, -434439541957702412239512$,\\
&\qquad$-151897369076330346742203, 610281064922154821463337, 211710160779065690463654, -678145211065889672618174$,\\
&\qquad$-230622710469268836012522, 587757024818899676327921, 191925680146871604212121, -391624630310001080455072$,\\
&\qquad$-119033636843278373618263, 197766785967009160429706, 53612258357612675910509, -74608688043563231273283$,\\
&\qquad$-17058617253176590894815, 20687510333056339496635, 3713992963764087642008, -4130593027985267121538$,\\
&\qquad$-531046062111838509185, 578162928492769417198, 46961262145834098724, -54781357679840005659, -2287539667564813302$,\\
&\qquad$\phantom{-}3352222063105017682, 40620243008078108, -123514842854460529, 1010652258526922, 2420447510811051, -55851319642213$,\\
&\qquad$-18452351744696, 786556882087$\\
\hline
\end{tabular}
}
\end{table}

\setlength{\tabcolsep}{7pt} 
\begin{table}[!tbh]
\caption[]{Continuation of Table 4.}
 \label{tab:tab05}
{\tiny
  \begin{tabular}{c l}
\hline  
$p,\Delta_K$  & Polynomial coefficients for decreasing powers of $x$\\
\hline
$73,293^{72}$&  $1, 1, -144, -139, 9730, 9056, -410441, -367906, 12132108, 10460357, -267391523, -221468185, 4565619563, 3627532367$,\\
(178) &\qquad$-61950350404, -47147580804, 679947392088, 494890787591, -6113510799041, -4248263015489, 45441140233604$,\\
&\qquad$\phantom{-}30093580771161, -281045946680083, -177038034694018, 1452868601357593, 868715820105921, -6295670418447183$,\\
&\qquad$-3565186787476080, 22902199637195106, 12253312487223193, -69962703250113206, -35272462117714863$,\\
&\qquad$\phantom{-}179338461559168548, 84955105196536588, -385048845118205565, -170849729564878739, 690518720387088459$,\\
&\qquad$\phantom{-}285992431892776900, -1030365887765938986, -396815828472409411, 1273033339340609215, 453948542708929508$,\\
&\qquad$-1294487651212092317, -425372577783968644, 1075460998363558124, 323938213588417248, -723673836702904009$,\\
&\qquad$-198620283444369134, 390347789427434985, 96972004303780760, -166728320860948405, -37205870049815405$,\\
&\qquad$\phantom{-}55579282507755157, 11040773363511690, -14210001829233769, -2484259728799330, 2727333323515561$,\\
&\qquad$\phantom{-}413148457058309, -382313978331429, -49086744961361, 37718946870926, 3978758518003, -2483523023306$,\\
&\qquad$-206824721197, 100506678700, 6386688064, -2177670239, -106831868, 20248339, 685305, -67318, -1667, 55, 1$\\
$79,317^{78}$&  $1, 1, -156, -151, 11476, 10742, -529833, -479120, 17237859, 15043356, -420695673, -353920938, 8005510587, 6484821551$,\\
(196)&\qquad$-121897544509, -94959796018, 1512506583433, 1131648317743, -15496341553928, -11120185043808, 132376432838106$,\\
&\qquad$\phantom{-}90975356804502, -949633381885839, -624050367056385, 5750894189811100, 3607649713387208, -29508660802301222$,\\
&\qquad$-17639571023194681, 128602029542560426, 73114526892928284, -476661570350174302, -257209659492690357$,\\
&\qquad$\phantom{-}1503149344492126347, 768122531281080170, -4030993135188740266, -1945953639652130387, 9180641646273875633$,\\
&\qquad$\phantom{-}4175544021158828515, -17720410438845095950, -7570210080186966239, 28904111895399788759, 11556972210324882738$,\\
&\qquad$-39694010615804869681, -14790824727072361901, 45687502119081734033, 15779507887426374013, -43834207774512448337$,\\
&\qquad$-13933215738460480507, 34832294570965230505, 10092292269108676458, -22752966081378874965, -5930089941589565108$,\\
&\qquad$\phantom{-}12110606741300452292, 2787350977545282191, -5198631736503715705, -1029729865708380737, 1777818076567389569$,\\
&\qquad$\phantom{-}292332297736928516, -477200723807927472, -61925190559965783, 98688594363464519, 9404904692311374$,\\
&\qquad$-15354866035626205, -967602418893197, 1742477355890453, 62117408755516, -138523447421396, -2288943125790$,\\
&\qquad$\phantom{-}7331249803128, 67426164550, -243056280171, -3889298793, 4696185844, 157596564, -45517649, -2527289, 140716, 12298, 218, 1$\\
$83,167^{82}$&  $1, 1, -82, -81, 3240, 3160, -82160, -79079, 1502501, 1426425, -21111090, -19757815, 237093780, 218618940, -2186189400$,\\
(183)&\qquad$-1984829850, 16871053725, 15071474661, -110524147514, -97082021465, 621324937376, 536211932256, -3022285436352$,\\
&\qquad$-2560547383576, 12802736917880, 10638894058520, -47465835030320, -38650751381832, 154603005527328, 123234279768160$,\\
&\qquad$-443643407165376, -345780890878896, 1123787895356412, 855420636763836, -2515943049305400, -1867897112363100$,\\
&\qquad$\phantom{-}4981058966301600, 3601688791018080, -8719878125622720, -6131164307078475, 13488561475572645, 9206478467454345$,\\
&\qquad$-18412956934908690, -12176310231149295, 22138745874816900, 14154280149473100, -23385332420868600$,\\
&\qquad$-14420954992868970, 21631432489303455, 12832205713993575, -17451799771031262, -9929472283517787$,\\
&\qquad$\phantom{-}12220888964329584, 6646448384109072, -7384942649010080, -3824345300380220, 3824345300380220, 1877405874732108$,\\
&\qquad$-1683191473897752, -779255311989700, 623404249591760, 270533919634160, -191991813933920, -77535155627160$,\\
&\qquad$\phantom{-}48459472266975, 18053528883775, -9847379391150, -3348108992991, 1575580702584, 482320623240, -192928249296$,\\
&\qquad$-52251400851, 17417133617, 4076350421, -1101716330, -215553195, 45379620, 7059052, -1086008, -123410, 12341, 861, -42, -1$\\
$89,179^{88}$&  $1, 1, -88, -87, 3741, 3655, -102340, -98770, 2024785, 1929501, -30872016, -29034396, 377447148, 350161812, -3801756816$,\\
(199)&\qquad$-3477216600, 32164253550, 28987537150, -231900297200, -205811513765, 1440680596355, 1258315963905, -7778680504140$,\\
&\qquad$-6681687099710, 36749279048405, 31022118677225, -152724276564800, -126600387152400, 560658857389200, 456002537343216$,\\
&\qquad$-1824010149372864, -1454278362337824, 5271759063474612, 4116305022165108, -13559593014190944, -10358022441395860$,\\
&\qquad$\phantom{-}31074067324187580, 23196134763125940, -63484158299081520, -46252743903616536, 115631859759041340, 82115378669464140$,\\
&\qquad$-187692294101632320, -129728497393775280, 271250494550621040, 182183167981760400, -348524321356411200$,\\
&\qquad$-227068876035237600, 397370533061665800, 250649105469666120, -401038568751465792, -244382877832924467$,\\
&\qquad$\phantom{-}357174975294274221, 209769429934732479, -279692573246309972, -157890968768078210, 191724747789809255$,\\
&\qquad$\phantom{-}103719945525634515, -114449595062769120, -59132290782430712, 59132290782430712, 29065024282889672$,\\
&\qquad$-26252279997448736, -12220888964329584, 9929472283517787, 4355031703297275, -3167295784216200$,\\
&\qquad$-1300853625660225, 841728816603675, 321387366339585, -183649923622620, -64617565719070, 32308782859535$,\\
&\qquad$\phantom{-}10363194502115, -4481381406320, -1292706174900, 476260169700, 121399651100, -37353738800, -8217822536$,\\
&\qquad$\phantom{-}2054455634, 377348994, -73629072, -10737573, 1533939, 163185, -15180, -990, 45, 1$\\
$97,389^{96}$&  $1, 1, -192, -187, 17578, 16664, -1021929, -942362, 42389952, 37997337, -1336462731, -1163680149, 33314897667, 28156698311$,\\
(249)&\qquad$-674408368288, -552838527076, 11301013598482, 8977909375093, -158994886726653, -122308729912222, 1898473578664977$,\\
&\qquad$\phantom{-}1412901222447415, -19399951044988978, -13955275813688939, 170762291665699972, 118614446384719822$,\\
&\qquad$-1301309882049061094, -871945780122664925, 8619377666392013922, 5565181768869578253, -49771732711929844788$,\\
&\qquad$-30930479393262936612, 251112071197806320739, 150020160792137684755, -1108677862726332708568$,\\
&\qquad$-635934910830427844783, 4287605330953710413992, 2358093640258069799412, -14530612000190559702170$,\\
&\qquad$-7651494316042985990505, 43150759530518584124075, 21722212320251222385725, -112228827609141155647105$,\\
&\qquad$-53922076560098598529529, 255396998183926539104038, 116914878588220059051162, -507829409247920514083495$,\\
&\qquad$-221083325574365962264596, 880677517892705249616378, 363892311994091133884903, -1329003780593599646666617$,\\
&\qquad$-520077115693458007818368, 1740436761368591469099430, 643542366602780582216934, -1971574632911812727476760$,\\
&\qquad$-687085946667407137439777, 1924701872709234943427512, 630420911347065941491567, -1612252793611361847320163$,\\
&\qquad$-494786845657867267549369, 1153114327395849183265703, 330399291621281386173050, -700204009654308355534634$,\\
&\qquad$-186550016549016488099855, 358661632153989150215180, 88424127745933995583848, -153831904737635556704045$,\\
&\qquad$-34895855606206393994169, 54780945797382119086423, 11357456950456450852662, -16039377858795390694256$,\\
&\qquad$-3015740835790853975136, 3817489549458070523980, 645389203032607765111, -728761616821001045732$,\\
&\qquad$-109825647572795780485, 109819500917589950691, 14644442153531449279, -12813785741765692684$,\\
&\qquad$-1506241051689218153, 1130523247638411169, 117451448799052895, -73215396010936265, -6798178894924984$,\\
&\qquad$\phantom{-}3351149496877315, 283215594877129, -103171823466943, -8071372448439, 1999976788345, 145004553360$,\\
&\qquad$-22289089940, -1467714628, 122419905, 7058361, -228410, -10990, 73, 1$\\
\hline
\end{tabular}
}
\end{table}

\end{document}